




\input amstex
\magnification=1200
\documentstyle{amsppt}
\baselineskip 20pt

\hsize 6.25truein
\vsize = 8.75truein
\def\cplus{\hbox{$\subset${\raise1.05pt\hbox{\kern -0.55em
${\scriptscriptstyle +}$}}\ }}

\def\bcplus{\hbox{$\supset${\raise1.05pt\hbox{\kern -0.55em
${\scriptscriptstyle +}$}}\ }}

\catcode`\@=11
\def\logo@{}
\catcode`\@=13

\centerline{\bf Generalized Harish-Chandra Modules: A New Direction}
\centerline{\bf in the Structure Theory of Representations}
\vskip .25in
\centerline{\bf Ivan Penkov, Gregg Zuckerman}

\vskip .20in \footnote""{Math. Subject Classification 2000:
Primary 17B10, Secondary 22E46, 22E47}

\document
\widestnumber\key{DMPP}
\def\supp{\text{\rm supp}}
\def\C{\bold C}
\def\R{\bold R}

\def\Z{\bold Z}

\vskip .20in \topmatter \abstract  Let $\frak g$ be a reductive
Lie algebra over $\bold C$.  We say that a $\frak g$-module $M$
is a generalized Harish-Chandra module if, for some subalgebra
$\frak k \subset\frak g$, $M$ is locally $\frak k$-finite and has
finite $\frak k$-multiplicities.  We believe that the problem of
classifying all irreducible generalized Harish-Chandra modules
could be tractable.  In this paper, we review the recent success
with the case when $\frak k$ is a Cartan subalgebra.  We also
review the recent determination of which reductive in $\frak g$
subalgebras $\frak k$ are essential to a classification. Finally,
we present in detail the emerging picture for the case when
$\frak k$ is a principal $3$-dimensional subalgebra.
\endabstract
\endtopmatter

\centerline{\bf INTRODUCTION}
\vskip .15in

Deep mathematical theories are usually rooted in a combination of
ideas.  It is also true that the core of mathematics, whatever it
may be, consists of theories which, due to their complexity, have
taken a long time to mature.  Representation theory is a perfect
illustration of both of these statements. By the end of the
$20^{th}$ century, representation theory has grown to an enormous
subject, with many different aspects and with complex relations
to theoretical physics, and with flavors ranging from
combinatorics, through abstract algebra, algebraic geometry,
homological algebra, to harmonic analysis and mathematical
physics.

A central part of the foundation of representation theory is the
Cartan-Killing classification of finite dimensional complex
simple Lie  algebras, or equivalently of all connected reduced
Dynkin diagrams.  Moreover, many other fundamental results which
have shaped the face of representation theory are also
classifications.  Strictly speaking, the representation theory of
Lie groups or Lie algebras starts with Cartan's classification of
irreducible finite dimensional representations, i.e. with the
classification of integral dominant weights.  One may say that
the skeleton of representation theory consists of several
explicit results and classifications such as H. Weyl's character
formula, the formula for the multiplicity of one Verma module in
another (commonly referred to as the Kazhdan-Lusztig conjecture)
the classification of Harish-Chandra modules, the classification
of simple Lie superalgebras, the constructions of Kac-Moody
algebras and quantum groups, etc.  Some classification problems
have been agreed upon to be unrealistic: a classical example is
the problem of classifying all irreducible representations of a
simple Lie algebra of rank greater than $1$.

We believe that in the last decade a combination of conceptual
developments has led to a possibility to restate this latter
problem, in a more restrictive but still enormously general way,
and to turn it into a tractable problem.  More precisely, let
$\frak g$ be a fixed complex reductive Lie algebra.  We are
interested in the problem of classifying all irreducible $\frak
g$-modules $M$ which have finite multiplicities as $\frak g
[M]$-modules, where $\frak g [M] \subset \frak g$ is the
subalgebra of all elements of $\frak g$ acting locally finitely
on $M$, i.e. $g\in\frak g [M]$, if, for any $m\in M$, the span
$\langle m,g\cdot m, g^2\cdot m,\ldots \rangle_{\C}$ is a finite
dimensional subspace of $M$.

The purpose of the present paper is to provide a brief review of
interrelated concepts and results which have led to this problem,
to describe the status quo, and to present some recent new
results.  Clearly, Harish-Chandra modules play a central role in
the subject, but we do not present them in this paper.  Two
fundamental references on Harish-Chandra modules are  \cite{V}
and \cite{KV}. We also omit all proofs given elsewhere in the
literature.

\vskip .20in

{\bf Acknowledgement}
\vskip .15in

We thank Professors G. van Dijk and V. Molchanov for giving us
the opportunity to present our results at the conference
``Representations of Lie Groups, Harmonic Analysis on Homogeneous
Spaces and Quantization'' in Leiden, December 2002.  We are also
grateful to Professor V. Serganova for some very interesting and
stimulating discussions of the subject. Finally, we thank
Professor J. Willenbring for the permission to state his
unpublished joint result with the second named author. \vskip
.15in

{\bf Notational conventions}
\vskip .15in

The ground field is $\bold C$.  We set $\bold Z_+ := \{n\geq 0
\mid n\in\bold Z\}$.  By $\langle\quad\rangle_{\bold Z_{+}}$,
$\langle\quad\rangle_{\bold R_{+}}$ and
$\langle\quad\rangle_{\bold C}$ we denote linear span over $\bold
Z_+, \bold R_+$, or $\bold C$.  The superscript $*$ indicates
dual space.  If $\frak l$ is a Lie algebra, $Z(\frak l)$ stands
for the center of $\frak l$, and if $\frak l$ is a Lie subalgebra
in a fixed Lie algebra $\frak g$, $C(\frak l)$ and $N(\frak l)$
denote respectively the centralizer and the normalizer of $\frak
l$ in $\frak g$. If $\frak l$ is reductive in $\frak g, C(\frak
l)$ is also reductive in $\frak g$.  The signs $\cplus$ and
$\bcplus$ stand for semidirect sum of Lie algebras.  If $\frak l$
is a reductive Lie algebra, we set $\frak l_{ss}: = [\frak
l,\frak l]$. \vskip .20in

\centerline{\bf 1.  General discussion and statement of the problem}
\vskip .10in

Let $\frak g$ be a reductive Lie algebra.

\proclaim{Theorem 1} Let $M$ be any $\frak g$-module.  Then the
set \newline $\frak g [M] := \{ g\in\frak g\mid g ~ \text{acts
locally finitely on} ~ M\}$ is a Lie subalgebra of $\frak g$.
\endproclaim

This result has been proved independently and by different
methods by V. Kac in \cite{K} and by S. Fernando in \cite{F}.  We
call $\frak g [M]$ the {\it Fernando-Kac} subalgebra of $M$, cf.
\cite{PSZ}.

Let $\frak k \subset\frak g$ be a Lie subalgebra.  We define a
$(\frak g,\frak k)$-{\it  module} $M$ as a $\frak g$-module $M$
such that $\frak k\subset\frak g [M]$.  A $(\frak g,\frak
k)$-module $M$ is {\it  strict} if $\frak k =\frak g [M]$.  If
$M$ is a $(\frak g,\frak k)$-module and $N$ is a finite
dimensional irreducible $\frak k$-module (or a $\frak k$-{\it
type} for short), we define the multiplicity $[M:N]$ as the
supremum of the multiplicities $[M':N]$ for all finite
dimensional $\frak k$-submodules $M'\subset N$.  We say that $M$
is of {\it  finite type over} $\frak k$ if $[M:N] < \infty$ for
any $N$, and we say that $M$ is of {\it  infinite type over}
$\frak k$ if $[M:N]\ne 0$ implies $[M:N] = \infty$ for all $N$.
In \cite{PS2} the following important technical result is proved.

\proclaim{Proposition 1}  Let $M$ be an irreducible $\frak
g$-module and $\frak k$ be a reductive in $\frak g$ subalgebra
with $\frak k \subset\frak g [M]$.  Then $\frak k$ acts
semisimply on $M$, and $M$ has either finite or infinite type
over $\frak k$.  Therefore, there is a canonical isomorphism of
$\frak k$-modules
$$
M\cong \bigoplus\limits_{W\in\hat \frak k} \text{Hom}_{\frak k}
(W,M)\otimes W,
$$
where $\hat{\frak k}$ is the set of equivalence classes of $\frak
k$-types, and the isotypic components \newline $\text{Hom}_{\frak
k}(W,M)\otimes W$ are either all finite dimensional or all
infinite dimensional.
 \endproclaim

We say that a subalgebra $\frak l\subset\frak g$ is of {\it
finite type} if there exists an irreducible strict $(\frak
g,\frak l)$-module $M$ of finite type over $\frak l$.

A classical theorem of Harish-Chandra implies that if a
subalgebra $\frak k$ equals the fixed points of an involution on
$\frak g$ (in this, the so-called symmetric case, $\frak k$ is
necessarily reductive in $\frak g$), then any irreducible $(\frak
g,\frak k)$-module has finite type over $\frak k$.  Assume $\frak
g$ is simple and $\frak k$ is symmetric.  In this case the only
possible subalgebras $\frak l$ with $\frak k\subset\frak
l\subset\frak g$ are parabolic, and it is easy to see that an
irreducible $(\frak g,\frak k)$-module $M$ is either a highest
weight module or a strict $(\frak g,\frak k)$-module, i.e. in the
latter case $\frak g [M] = \frak k$. For the purpose of this
paper, we define a {\it  Harish-Chandra module} as a $(\frak
g,\frak k)$-module of finite type over a symmetric subalgebra
$\frak k\subset \frak g$.

Classically Harish-Chandra modules are defined as $(\frak g,
K)$-modules, where $K$ is a subgroup of an algebraic group $G$
and the Lie algebra $\frak k$ of $K$ is a symmetric subalgebra of
the Lie algebra $\frak g$ of $G$. This definition is closely
related to the above definition, but is not equivalent to it.  In
the original definition, Harish-Chandra modules have been
classified by a monumental effort of several groups of
mathematicians including R. Langlands, D. Vogan, A. Beilinson, J.
Bernstein and others.  The literature on the subject is enormous;
see the texts \cite{V} and \cite{KV} (and the references therein)
for a presentation and discussion of the classification.

Harish-Chandra modules have their roots in physics and have been
recognized to be of fundamental importance because of their
specific properties, in particular their relation to
unitarizability of representations, and not because of a clear
understanding of the place of Harish-Chandra modules among
general $\frak g$-modules.  In our opinion, such an understanding
can be based on the notion of the Fernando-Kac subalgebra, which
was not part of the original theory.  Moreover, this notion
enables us to consider Harish-Chandra modules and weight modules,
two seemingly unrelated subjects, from a single point of view.

We define a {\it  generalized Harish-Chandra} $\frak g${\it
-module} to be a $\frak g$-module $M$ which has finite type over
a subalgebra of $\frak g [M]$.  An example of a generalized
Harish-Chandra module is a weight module (the definition is
recalled in subsection 2.1 below) with finite weight
multiplicities.  The theory of weight modules has developed
practically independently with that of Harish-Chandra modules and
has culminated in O. Mathieu's classification of irreducible
weight modules with finite weight multiplicities, \cite{M}.  We
present a summary of Mathieu's result in subsection 2.2.

We believe that the classification of Harish-Chandra modules,
together with Mathieu's classification suggest that the problem
of classifying all irreducible generalized Harish-Chandra modules
(posed in slightly different terms in \cite{PSZ}) could be
tractable.  If that is the case, this classification problem
could be the ultimate substitute for the unrealistic problem of
classifying all irreducible $\frak g$-modules.  Closely related
problems are as follows. \vskip .10in

{\bf A.}  For a given subalgebra $\frak l\subset\frak g$ of
finite type, classify all irreducible $(\frak g,\frak l)$-modules
of finite type over $\frak l$. \vskip .10in

{\bf B.}  For a given subalgebra $\frak k\subset \frak g$,
reductive in $\frak g$ and of finite type, classify all
irreducible
 $(\frak g,\frak k)$-modules $M$ of finite type over $\frak k$ with $\frak g [M]_{\text{red}} = \frak k$.
\vskip .10in

{\bf C.}  For a given Fernando-Kac subalgebra $\frak
l\subset\frak g$ of finite type, classify all irreducible strict
$(\frak g, \frak l)$-modules of finite type over $\frak g$.
\vskip .10in

No systematic solution of any of the above problems is known
except in the cases of Harish-Chandra modules and weight modules.
In this paper we give a brief account of some of the developments
which have led naturally to those problems.  In section 2 below,
we present a summary of results on weight modules.  In section 3,
we summarize some recent results of \cite{PS2} and \cite{PSZ} on
the description of Fernando-Kac subalgebras of finite type and
related subjects.  In particular we give an explicit description
of any reductive in $\frak g$ subalgebra $\frak k$ which is the
reductive part of a Fernando-Kac subalgebra $\frak l\subset\frak
g$ of finite type.

In section 4, we discuss the case when $\frak g$ is simple and
$\frak k = \frak{sl} (2)$ is a principal $\frak{sl}
(2)$-subalgebra.  As the principal $\frak{sl} (2)$-subalgebra is
never symmetric, unless $\frak g\cong \frak{sl}(2)$ or
$\frak{sl}(3)$, the corresponding $(\frak g,\frak k)$-modules of
finite type over $\frak k$ are interesting objects of study in
the new theory of generalized Harish-Chandra modules.  Their
existence is ensured by a result of \cite{PSZ}, see Theorem 4
below. A careful study of this result in the case when $\frak k$
is a principal $\frak{sl} (2)$-subalgebra, shows that the
irreducible strict $(\frak g,\frak k)$-modules $M$ constructed in
\cite{PSZ} are in some sense generic, and not every given $\frak
k$-type occurs in some  module $M$.

In the present paper, we prove two new results: that any infinite
dimensional irreducible $(\frak g,\frak k)$-module of finite type
over $\frak k$ is strict, and that any fixed $\frak k$-type $Z$
occurs in some irreducible $(\frak g,\frak k)$-module as the
$\frak k$-type of minimal dimension.  The second result suggests
that some aspects of  Vogan's theory of minimal $\frak k$-types
may carry over to generalized Harish-Chandra modules. 

\newpage

\centerline{\bf 2.  The case when $\frak h\subset\frak g [M]$:
irreducible weight modules} \vskip .10in

The case when $\frak h\subset\frak g [M]$ is an important case in
which much more is known than in the general case.  Below we
present a very brief survey of known results. \vskip  .20in

{\bf 2.1. General invariants.} In this subsection $M$ is an
arbitrary irreducible $\frak g$-module such that $\frak
h\subset\frak g [M]$.  We do not assume that $M$ is of finite
type over $\frak g [M]$ or over $\frak h$. Let $\Delta$ denote
the roots of $\frak g$ with respect to $\frak h$, and
$\Delta_{\frak g [M]}$ denote the roots of $\frak g [M]$ with
respect to $\frak h$.

Note that if $\frak h\subset\frak g [M]$, then $\frak g [M]$ has
a unique reductive part $\frak g [M]_{\text{red}}$ (which is the
unique maximal reductive subalgebra of $\frak g [M]$ and is
automatically reductive in $\frak g$), and by Proposition 1,
$\frak g [M]_{\text{red}}$ acts semisimply on $M$.  Consequently,
$\frak h$ acts also semisimply on $M$ and as an $\frak h$-module
$M$ has the decomposition 
$$
M = \bigoplus\limits_{\nu\in\frak h^{*}} M^\nu, \tag1
$$
where the weight spaces $M^\nu := \{ m\in M\mid h\cdot m = \nu
(h)\cdot m, \forall h\in\frak h\}$ are all either finite
dimensional or infinite dimensional.  A $\frak g$-module which
admits a decomposition (1) is by definition a {\it  weight $\frak
g$-module}.  We set $\supp M:= \{ \nu\in\frak h^*\mid M^\nu\ne
0\}$. Let $\Gamma_M$ denote the submonoid of
$\langle\Delta\rangle_{\Z}$ generated by
$\Delta\backslash\Delta_{\frak g [M]}$.  We define the $M$-{\it
decomposition} of $\Delta$, or {\it  the shadow decomposition of
$\Delta$ corresponding to }$M$, to be the decomposition
$$
\Delta = \Delta^I_M\cup \Delta^F_M\cup \Delta^+_M\cup \Delta^-
_M, \tag2
$$
where
$$\align
&\Delta^I_M := \{ \alpha\in\Delta \mid \alpha\in \langle\Gamma_M\rangle_{\R_{+}}, - \alpha\in\langle\Gamma_M\rangle_{\R_{+}}\},  \\
&\Delta^F_M := \{ \alpha\in\Delta \mid \alpha\notin \langle\Gamma_M\rangle_{\R_{+}}, - \alpha\notin\langle\Gamma_M\rangle_{\R_{+}}\},  \\
&\Delta^+_M := \{ \alpha\in\Delta \mid \alpha\notin\langle\Gamma_M\rangle_{\R_{+}}, - \alpha\in\langle\Gamma_M\rangle_{\R_{+}}\},  \\
&\Delta^- _M := \{ \alpha\in\Delta \mid \alpha\in \langle\Gamma_M\rangle_{\R_{+}}, - \alpha\notin\langle\Gamma_M\rangle_{\R_{+}}\}.
\endalign
$$
In particular, the $M$-decomposition of $\Delta$ is determined
only by $\frak g [M]$.  The decomposition (2) induces a
decomposition of $\frak g$,
$$
\frak g = (\frak g^I_M + \frak g^F_M) \oplus \frak g^+_M\oplus
\frak g^- _M,
$$
where
$$
\frak g^I_M := \frak h \oplus
(\bigoplus\limits_{\alpha\in\Delta^{I}_{M}} \frak g^\alpha),\quad
\frak g^F_M := \frak h\oplus
(\bigoplus\limits_{\alpha\in\Delta^{F}_{M}} \frak g^\alpha),
$$
$$
\frak g^\pm _M := \bigoplus\limits_{\alpha\in\Delta^{\pm}_{M}}
\frak g^\alpha .
$$
It follows from the main result of \cite{DMP} that $\frak p_M :=
(\frak g^I_M + \frak g^F_M)\oplus\frak g^+_M$ is a parabolic
subalgebra whose semisimple part is the direct sum $[\frak g^I_M,
\frak g^I_M] \oplus [\frak g^F_M, \frak g^F_M]$.

Note that the shadow decomposition of $\Delta$ corresponding to
$M$ depends only on $\frak g [M]$.  In general it is not true
that the shadow decomposition reconstructs $\frak g [M]$.  One
can only show that
$$
\frak g [M] = (\frak g^F_M + (\frak g^I_M\cap \frak g [M]))
\oplus \frak g^+ _M,
$$
and the results of \cite{PS2} and \cite{PSZ} imply that there are
irreducible weight modules with the same shadow decomposition and
different Fernando-Kac subalgebras, see subsection  3.3 below.

The shadow decomposition reconstructs ``the shape'' of $\supp M$,
more precisely it reconstructs $\supp ~ M$ up to adding the
support of an arbitrary irreducible (finite dimensional) $\frak
g^F_M\oplus\frak g^+_M$-submodule of $M$.  Indeed, a direct
argument, see \cite{PS1}, shows that
$$
\supp M = \supp M^F + \Gamma_M,
$$
where $M^F$ is any irreducible $\frak g^F_M\oplus \frak
g^+_M$-submodule of $M$.

In general, irreducible weight modules $M$ are not generalized
Harish-Chandra modules (as they are not necessarily of finite
type over $\frak g [M]$) and no classification of irreducible
weight modules is available.  The following theorem provides a
reduction, which can be extended to a classification in the case
when $M$ has finite $\frak h$-type.

\proclaim{Theorem 2}  $M$ has a unique irreducible $\frak
p_M$-submodule $M^{IF}$ on which $\frak g^+_M$ acts trivially.
Furthermore, there is an isomoprhism of $(\frak g^I_M + \frak
g^F_M)$-modules $M^{IF}\simeq M^I\otimes M^F$, where $M^I$ is an
irreducible $\frak g^I_M$-module with $\frak g^I_M = (\frak
g^I_M)^I_{M^{I}}$ and $M^F$ is an irreducible finite dimensional
$\frak g^F_M$-module.  Finally, $M$ is the unique irreducible
quotient of the induced $\frak g$-module $U(\frak
g)\bigotimes\limits_{U(\frak p_{M})} M^{IF}$.
\endproclaim

The proof of Theorem 2 see in \cite{DMP}.

Theorem 2 reduces the study of an arbitrary irreducible $\frak
g$-module $M$ with $\frak h \subset \frak g [M]$ to the study of
the irreducible $\frak g^I_M$-module $M^I$.  The latter module is
a {\it  cuspidal} $\frak g^I_M$-module, which by definition means
that its shadow decomposition is trivial in the sense that
$\Delta^I = (\Delta^I)^I_{M^{I}}$. Cuspidal modules arise in
Harish-Chandra module theory, for if $\frak g$ is simple, any
Harish-Chandra module $M$ for which $\frak h\subset\frak g [M]$
and such that $M$ is not a highest (or a lowest) weight module,
is necessarily a cuspidal $\frak g$-module.

Mathieu's classification result, which we describe in the next
subsection, provides an explicit classification of all cuspidal
$\frak g$-modules with finite dimensional weight spaces (i.e. of
finite type over $\frak h$), and via Theorem 2 this yields a
complete classification of all irreducible generalized
Harish-Chandra modules $M$ such that $\frak h\subset \frak g [M]$
and $M$ has finite type over $\frak h$.  No classification of
irreducible generalized Harish-Chandra modules which have
infinite type over $\frak h$ is known, except when they are
Harish-Chandra modules. \vskip .20in

{\bf 2.2 Fernando's theorem and Mathieu's classification.}  In
\cite{F}, S. Fernando constructed the shadow decomposition (2) of
any irreducible weight module $M$ of finite type over $\frak h$,
proved Theorem 2 in that case, and showed that the Fernando-Kac
subalgebra of $M$ is determined by the shadow decomposition via the formula:
$$
\frak g [M] = \frak g^F_M \oplus \frak g^+ _M.
$$
In particular, if $M$ is cuspidal, and of finite type over $\frak
h$, we have $\frak g [M] = \frak h$.  Moreover, in this case all
weight spaces $M^\nu$ are immediately seen to be of the same
dimension $d$; by definition, $M$ is called then a cuspidal
module of {\it  degree} $d$.  In \cite{F} Fernando established
also the following key result.

\proclaim{Theorem 3} The reductive Lie algebra $\frak g$ admits
an irreducible cuspidal weight module of finite type over $\frak
h$ if and only if all simple components of $\frak g_{ss}$ are of
type $A$ and $C$. \endproclaim

The further trivial observation that any irreducible cuspidal
module over $\frak g$ is isomorphic to a tensor product of
cuspidal irreducible modules over the simple components of $\frak
g_{ss}$ with a $1$-dimensional module over the center of $\frak
g$, reduces the problem of classifying all cuspidal irreducible
$\frak g$-modules of finite type over $\frak h$ to the same
problem for the Lie algebras $\frak{sl} (n+1)$ and
$\frak{sp}(2n)$.   This latter problem was solved completely in
Mathieu's breakthrough paper \cite{M}.

Mathieu's main idea is that irreducible cuspidal weight modules
come in coherent families, each family being determined by a
finite dimensional, irreducible representation of a corresponding
maximal reductive root subalgebra.  The classification then
reduces to describing a continuous parameter (the position of the
module within the family) and a mixed (partly continuous, partly
discrete) parameter (the highest weight of an irreducible
representation of a reductive Lie algebra).

In the rest of this section $\frak g = \frak{sl}(n+1)$,
$\frak{sp}(2n)$.  Following Mathieu we call a (reducible) weight
$\frak g$-module $\Cal M = \bigoplus\limits_{\nu\in\frak h^{*}}
\Cal M^\nu$ a {\it  coherent family of degree $d$} if $\supp ~
\Cal M = \frak h^*$, $\dim \Cal M^\nu = \dim \Cal M^\mu = d$, and
for any $u$ in the centralizer of $\frak h$ in $U(\frak g)$, the
function $\lambda\in\frak h^*\mapsto tr ~ u\mid_{\Cal
M^{\lambda}}$ is a polynomial in $\lambda$. $\Cal M$ is called
{\it  semisimple} if it is semisimple as a $\frak g$-module.

A {\it  fiber} of the family $\Cal M$ is a $\frak g$-submodule
$\Cal M'$ of the form
 $\bigoplus\limits_{\gamma\in\langle\Delta\rangle_{\Z}} \Cal M^{x
+ \gamma}$ for some $x\in \frak h^*$.  Clearly $\Cal M$ is the
direct sum of all its fibers.  By definition, $\Cal M$ is an {\it
irreducible coherent} family if at least one fiber of $\Cal M$ is
irreducible (then necessarily almost any fiber of $\Cal M$ is
irreducible, see Lemma 4.7 in \cite{M}).

Mathieu's result is that for each irreducible cuspidal weight
module $M$ of degree $d$, there is a unique (up to isomorphism)
semisimple coherent family $\tilde M$ of degree $d$ for which $M$
is a fiber of $\tilde M$.  Furthermore, $\tilde M$ has a
(non-unique) irreducible infinite dimensional highest weight
submodule $L(\lambda_M)$ (with respect to a fixed Borel
subalgebra of $\frak g$) which necessarily has the property that
the dimensions of its weight spaces are bounded.  More generally
we will say that a weight module $M$ has {\it  bounded
multiplicities} if, for some $\ell$, $\dim M^\nu < \ell$ for all
$\nu\in\supp ~M$.\footnote{Mathieu uses the term admissible
weight module but, since this term is a synonym for a $(\frak
g,\frak k)$-module of finite type in Harish-Chandra module
theory, we prefer the term weight module with bounded
multiplicities.}

 The classification is based on three main facts:
\roster \item"{}" - that $L(\lambda_M)$ determines $\tilde M$ up
to isomorphism; in what follows we will write also
$\widetilde{L(\lambda_M)}$ for $\tilde M$; \vskip .06in

\item"{}" - that the correspondence $L(x)\mapsto
\widetilde{L(x)}$ is defined for any infinite dimensional highest
weight module $L(x)$ with bounded multiplicities; \vskip .06in

\item"{}" - that one can describe explicitly all weights $x$ for
which $\dim L(x) = \infty$ and $L(x)$ has bounded multiplicities.
\endroster
\vskip .06in

Here is, for instance, an explicit description of all such
weights $x$ for $\frak g = \frak{sp} (2n)$. Fix a basis
$\varepsilon_1,\ldots ,\varepsilon_n$ of $\frak h^*$ such that
$\varepsilon_1-\varepsilon_2,\ldots,\varepsilon_{n-1}-\varepsilon_n,
2\varepsilon_n$ is a system of simple roots.

\proclaim{Proposition 2} (\cite{M}, Lemma 9.1) Let $\frak g =
\frak{sp} (2n)$.  If $x = \sum\limits_{i} x_i\varepsilon_i$ and
\newline $\dim L(x) = \infty$, then $L(x)$ has bounded
multiplicities if and only if $x_i\in \Z + \frac{1}{2}$ and
\newline  $x_1 > x_2 >\ldots > x_{n-1} > \vert x_n\vert$.
\endproclaim

For $\frak g = \frak{sp}(2n)$ all weights $x$ with the above
property form a discrete set.  For $\frak g = \frak{sl}(n+1)$
this is no longer true and the description is slightly more
complicated, see \cite{M}.

Mathieu's classification can now be stated from the opposite end
as follows. \vskip .10in

1.  Determine the highest weights $x$ of all infinite dimensional
irreducible modules $L(x)$ with bounded multiplicities.  Define
two such weights $x$ and $x'$ to be equivalent if
$\widetilde{L(x)}\simeq\widetilde{L(x')}$.  Then the set of
equivalence classes parametrizes all semisimple coherent families. In
particular, for $\frak g = \frak{sp}(2n)$ two weights $x = \Sigma
x_i\varepsilon_i$ and $x' = \Sigma x'_i\varepsilon_i$ as above
are equivalent if and only if $x'_n = \pm x_n$. \vskip .10in

2.  Let $\tilde x$ be the equivalence class of a weight $x$ as in
Step 1, and let $\Cal M = \widetilde{L(x)}$.  Describe the subset
in $\frak h^*/\langle\Delta\rangle_{\Z}$ for which the
corresponding fibers on $\Cal M$ are irreducible.  Mathieu gives
an explicit combinatorial description of this set in terms of
$\tilde x$ and shows in particular that its complement
(corresponding to reducible fibers) is always the union of
precisely $n+1$ codimension $1$ subsets in $\frak
h^*/\langle\Delta\rangle_{\Z}$ for $\frak g = \frak{sl}(n+1)$,
and respectively $n$ codimension $1$ subsets for $\frak g =
\frak{sp} (2n)$.  In the latter case, the corresponding
conditions on $\eta = (\eta_1,\ldots, \eta_n)\in\frak
h^*/\langle\Delta\rangle_{\Z}$ are simply $\eta_i\notin \Z +
\frac{1}{2}$. \vskip .10in

3.  Finally, in each case a weight $x$ as above determines a
natural maximal reductive root subalgebra $\frak g'_x\subset
\frak g$ for which the weight $x$ is dominant.  If $\frak g =
\frak{sl} (n+1)$, then $\frak g'_x$ is always isomorphic to
$\frak {gl}(n)$, while for $\frak g = \frak{sp}(2n)$, we have
$\frak g'_x = \frak o (2n)$.  Furthermore, the degree $d$ of the
coherent family $\widetilde{L(x)}$ is computed in terms of $x$
and $\frak g'_x$ as follows: \vskip .06in \roster \item"{}" - for
$\frak g = \frak{sp} (2n)$, $d = \frac{1}{2^{n-1}} \dim L_{\frak
g'_{x}} (x+\varepsilon)$, where $\varepsilon = \sum\limits_{i}
\varepsilon_i$; \vskip .06in \item"{}" - for $\frak g =
\frak{sl}(n+1)$, $d = \dim L_{\frak g'_{x}}(x)$ unless the
infinitesimal character of $L(x)$ is regular integral; in the
latter case $d$ is an alternating sum of dimensions of finite
dimensional $\frak g'_x$-modules.
\endroster
\vskip .06in

The details are in \cite{M}.

Mathieu's classification has been generalized by D. Grantcharov
\cite{G} to the case of the Lie superalgebra $\frak{sl}(n/1)$,
and by I. Dimitrov \cite{Di} to the case of the infinite
dimensional Lie algebra $\frak{gl}(\infty)$.

\vskip .15in

\centerline{\bf 3.  Fernando-Kac subalgebras of finite type}
\vskip .10in

A part of a future classification of irreducible generalized
Harish-Chandra modules should be a good understanding of their
respective Fernando-Kac subalgebras.  The problem of classifying
all Fernando-Kac subalgebras and all Fernando-Kac subalgebras of
finite type of a given reductive Lie algebra $\frak g$ has been
addressed in the recent papers \cite{PS2} and \cite{PSZ}.  Below
we present a summary of the results. \vskip .10in

{\bf 3.1.  General Fernando-Kac subalgebras.} Little is known
about a general description of Fernando-Kac subalgebras of
irreducible $\frak g$-modules which are not generalized
Harish-Chandra modules.   In \cite{PS2} the following two results
are established. \vskip .06in

1.  An example of a subalgebra of $\frak{sl}(n)$ which is not a
Fernando-Kac subalgebra is constructed. \vskip .06in

2.  It is proved that any subalgebra $\frak{l} \subset\frak g$
which contains a Cartan subalgebra of $\frak g$ is a Fernando-Kac
subalgebra.  The proof is an explicit $\Cal D$ -module
construction which provides a strict irreducible $(\frak g,\frak
l)$-module $L$. This result, together with Theorem 8 below,
implies that not every Fernando-Kac subalgebra is of finite type.
\vskip .10in

{\bf 3.2.  A description of primal subalgebras.} In \cite{PSZ} a
{\it  primal subalgebra} of $\frak g$ is defined as a reductive
in $\frak g$ subalgebra $\frak k$ for which there exists an
irreducible generalized Harish-Chandra module $M$ such that
$\frak k$ is a maximal reductive subalgebra of $\frak g [M]$.  A
main result of \cite{PSZ} is the following theorem. \vskip .10in

\proclaim{Theorem 4}  A reductive in $\frak g$ subalgebra $\frak
k$ is primal if and only if $\frak k\cap C(\frak k_{ss})$ is a
Cartan subalgebra of $C(\frak k_{ss})$, or equivalently if
$C(\frak k) = Z(\frak k)$.
\endproclaim
\vskip .10in

The proof of Theorem 4 is also an explicit $\Cal D$-module
construction which yields an irreducible generalized
Harish-Chandra module $M$ for which $\frak k$ is a maximal
reductive subalgebra in $\frak g [M]$; see Theorems 4.3 and 4.4
in \cite{PSZ}.

Theorem 4 gives an explicit description of all primal subalgebras
of $\frak g$.  Indeed, recall that all semisimple subalgebras of
a semisimple, or equivalently reductive, Lie algebra have been
classified in the fundamental papers \cite{D1} and \cite{D2} of
E. Dynkin.    Then, Theorem 4 implies that for any semisimple
subalgebra $\frak k' \subset \frak g$, the primal subalgebras
$\frak k$ with $\frak k_{ss} = \frak k'$ are precisely all direct
sums $\frak k'\oplus \frak h_{C(\frak k')}$, where $\frak
h_{C(\frak k')}$ is any Cartan subalgebra of $C(\frak k')$. In
particular, Theorem 4 implies that any semisimple subalgebra of
$\frak g$ is the semisimple part of a Fernando-Kac subalgebra of
finite type.  Another corollary of Theorem 4 is that every
maximal proper subalgebra $\frak l\subset\frak g$ is a
Fernando-Kac subalgebra of finite type.  For the proof, as well
as for other related corollaries of Theorem 4, see \cite{PSZ}.
\vskip .10in

{\bf 3.3.  Fernando-Kac subalgebras of finite type.}  In general,
the problem of describing all Fernando-Kac subalgebras of finite
type of a given reductive (or simple) Lie algebra $\frak g$ is
open.  A complete description of all Fernando-Kac subalgebras of
$\frak g$ is known under various additional conditions.  The
following general theorem is proved in \cite{PSZ}. \vskip .10in

\proclaim{Theorem 5}  Let $\frak l\subset\frak g$ be a
Fernando-Kac subalgebra of finite type. \vskip .06in

1.  $N(\frak l) = \frak l$; hence $\frak l$ is an algebraic
subalgebra of $\frak g$. \vskip .06in

2.  There is a decomposition $\frak l = \frak n_{\frak l} \cplus
\frak l_{\text{red}}$, unique up to an inner automorphism of
$\frak l$, where $\frak l_{\text{red}}$ is a (maximal) subalgebra
of $\frak l$ reductive in $\frak g$. \vskip .06in

3.  Any irreducible $(\frak g,\frak l)$-module $M$ of finite type
over $\frak l$ has finite type over $\frak l_{\text{red}}$, and
$\frak l_{\text{red}}$ acts semisimply on $M$. \vskip .06in

4.  $C(\frak l_{\text{red}}) = Z(\frak l_{\text{red}})$ and
$Z(\frak l_{\text{red}})$ is a Cartan subalgebra of $C(\frak
l_{ss})$. \vskip .06in

5.  $\frak l\cap C (\frak l_{ss})$ is a solvable Fernando-Kac
subalgebra of finite type of $C(\frak l_{ss})$.
\endproclaim
\vskip .10in

The following theorem provides a complete description of all
solvable Fernando-Kac subalgebras of finite type. \vskip .10in

\proclaim{Theorem 6}  A solvable subalgebra $\frak s\subset\frak
g$ is a Fernando-Kac subalgebra of finite type if and only if
$\frak s$ contains a Cartan subalgebra of $\frak g$ and its
nilradical $\frak n_{\frak s}$ is the nilradical of a parabolic
subalgebra of $\frak g$ whose simple components are all of type
$A$ and $C$.
\endproclaim
\vskip .10in

For the proof see Proposition 3.2 in \cite{PSZ}.

For $\frak g = \frak{gl}(n)$, $\frak{sl}(n)$ the  following
theorem provides a complete description of all reductive
Fernando-Kac subalgebras of finite type. \vskip .10in

\proclaim{Theorem 7} (Theorem 5.1 in \cite{PSZ})  A reductive in
$\frak g = \frak{gl}(n), \frak{sl}(n)$ subalgebra $\frak k$ is a
Fernando-Kac subalgebra if and only if it is primal. (An explicit
description of primal subalgebras is provided by Theorem 4).
\endproclaim
\vskip .10in

For simple Lie algebras other than of type $A$, it is not known
whether every primal subalgebra is a Fernando-Kac subalgebra of
finite type.  The following proposition is proved in \cite{PS2}
and provides a partial answer to this question. \vskip .10in

\proclaim{Proposition 3}  Let $\frak g$ be simple of type other
than $B_n, n\geq 3$, and $F_4$.  Then any reductive root subalgebra (which is
automatically primal by Theorem 4) is a Fernando-Kac subalgebra
of finite type.
\endproclaim
\vskip .10in

We complete this section by a description of all Fernando-Kac
 subalgebras of finite type of $\frak g = \frak{gl}(n),
\frak{sl}(n)$ which are root subalgebras.  Let $\frak g =
\frak{gl}(n),\frak{sl}(n)$, $\frak h$ be a fixed Cartan
subalgebra, and $\frak l \supset \frak h$ be a root subalgebra of
$\frak g$.  Then $\frak l$ is determined by its subset of roots
$\Delta (\frak l)\subset\Delta$.  Let $\frak k = \frak
l_{\text{red}}$ and $\frak l = \frak k \bcplus \frak n$.  Fix an
arbitrary Borel subalgebra $\frak b\subset\frak g$ containing
$\frak h$ and let $N$ be any finite dimensional semisimple $\frak
k$-module.  Denote by $\Cal S_{\frak k}(N)$ the set of weights of
all $\frak k\cap\frak b$-singular vectors in $N$, and put $\Cal
C_{\frak k}(N) := \langle\Cal S_{\frak k}(N)\rangle_{\bold Z_+}$.
 The following theorem is proved in \cite{PSZ} (Theorem 5.8).
\vskip .10in

\proclaim{Theorem 8}  A root subalgebra $\frak l = \frak k\bcplus
\frak n = \frak h$ of $\frak g = \frak{gl}(n), \frak{sl}(n)$ is a
Fernando-Kac subalgebra of finite type if and only if $\Cal
C_{\frak k}(\frak g/\frak l)\cap\Cal C_{\frak k} (\frak n) =
\{0\}$.
\endproclaim
\vskip .10in

It is impossible not to notice that the condition in Theorem 8 is
a generalization of a parabolic or ``triangular'' decomposition:
if $\frak k = \frak h$, it means precisely that $\frak n$ is the complement of
a parabolic subalgebra, cf. Theorem 6 above.  It would be very interesting to find the analog
of Theorem 8 for a general simple Lie algebra.
 \vskip .10in

Finally, note that  Theorem 8 implies that any reductive root
subalgebra $\frak k$ of $\frak g = \frak{gl}(n), \frak{sl}(n)$ is
a Fernando-Kac subalgebra of finite type. Furthermore, it is easy
to check that any strict $(\frak g,\frak k)$-module $M$ is
cuspidal and is necessarily of infinite type over $\frak h$
unless $\frak k=\frak h$.  In particular this shows, as claimed
in subsection 2.1 above, that in general the shadow decomposition
of an irreducible $\frak g$-module $M$ does not determine the
subalgebra $\frak g[M]\subset\frak g$. \vskip .15in

\centerline{\bf 4. A case when $\frak h\not\subset\frak g [M]$: a
principal $\frak{sl}(2)$-subalgebra} \vskip .06in

In this section, we consider a specific, yet broad enough, class
of generalized Harish-Chandra modules which has not been
discussed in the literature.

According to Theorem 4, if $\frak g$ is simple and $\frak
k\subset\frak g$ is a simple subalgebra with $C(\frak k) = 0$,
then $\frak k$ is primal.  If, in addition, we know that any
intermediate subalgebra $\frak l, \frak k\subset\frak
l\subset\frak g$, is reductive, Theorem 4 implies that $\frak k$
is itself a Fernando-Kac subalgebra of finite type.  An important
example of such a situation is when $\frak k$ is a principal
$3$-dimensional subalgebra.  More precisely, for any simple
$\frak g$, an injective homomorphism $\frak{sl}(2)\to\frak g$ is
called {\it principal} if its image contains a regular nilpotent
element.  Principal $\frak{sl}(2)$-subalgebras of simple Lie
algebras have been studied in detail by \cite{D2} and \cite{Ko}.
In particular, it is true that every intermediate subalgebra
$\frak k\subset\frak l\subset\frak g$ (where $\frak k$ is a
principal $\frak{sl}(2)$-subalgebra) is semisimple, see
Proposition 5 below.  Hence $\frak k$ is a Fernando-Kac
subalgebra of finite type, and certain irreducible strict $(\frak
g, \frak k)$-modules of finite type have been constructed in
\cite{PSZ}. In the present paper, we prove the following more
detailed result.

\proclaim{Theorem 9}  Let $\frak g$ be simple, of rank greater or
equal to $2$, and $\frak k\subset\frak g$ be a principal
$\frak{sl}(2)$-subalgebra. \roster \item"{(i)}" Any infinite
dimensional irreducible $(\frak g, \frak k)$-module of finite
type over $\frak k$ is strict. \vskip .06in

\item"{(ii)}"  Let $Z$ be a fixed $\frak k$-type.  There exists
an irreducible infinite dimensional (and thus strict) $(\frak
g,\frak k)$-module $X$ of finite type over $\frak k$ such that
$Z$ occurs in $X$ with multiplicity $1$ and for any other $\frak
k$-type $W$ which occurs in $X$, $\dim W > \dim Z$.
\endroster
\endproclaim

In the next subsection we prove (i).  To prove (ii) we need to
recall some basic facts about a cohomological method of
constructing $(\frak g,\frak k)$-modules.  We do this in
subsection 4.2.  Finally, we present the proof of (ii) in
subsection 4.3.

\vskip .10in

{\bf 4.1.  A result about finite dimensional modules and its
application to generalized Harish-Chandra modules. } J.
Willenbring and the second named author recently proved the
following fact.

\proclaim{Proposition 4}  Let $\frak k$ be an
$\frak{sl}(2)$-subalgebra of a semisimple Lie algebra $\frak s$,
none of whose simple factors is isomorphic to $\frak{sl} (2)$.
Then there exists a positive integer $b(\frak s)$, such that for
every irreducible finite dimensional $\frak s$-module $V$, there
exists an injection of $\frak k$-modules $W\to V$, where $W$ is
an irreducible $\frak k$-module of dimension less than $b(\frak
s)$.
\endproclaim

A more general statement will be proved in \cite{WZ}.  The proof
uses invariant theory.

To relate Proposition 4 with the statement of Theorem 9 (i), we
recall the notion of a stem subalgebra.  In a slight modification
of the original terminology of \cite{D2}, a {\it stem subalgebra}
is defined in \cite{PSZ} as a subalgebra $\frak s$ of a
semisimple Lie algebra $\frak r$ which is not contained in a root
subalgebra of $\frak r$.  The following result is established in
\cite{D2} (see p. 160).

\proclaim{Proposition 5}  Let $\frak s\subset\frak r$ be a stem
subalgebra.  Then any intermediate subalgebra $\frak s'$, $\frak
s\subset \frak s'\subset \frak r$, is a semisimple stem
subalgebra.  If $\frak r$ is simple and $\frak l$ is a principal $\frak{sl}(2)$-subalgebra
of $\frak r$, then $\frak k$ is a stem subalgebra of $\frak r$, and for any intermediate subalgebra $\frak s'$, no simple factor of $\frak s'$ is isomorphic to $\frak{sl}(2)$.
\endproclaim
\vskip .06in

Propositions 4 and 5 imply Theorem 9 (i).  Indeed, consider any
infinite dimensional irreducible $(\frak g,\frak k)$-module $M$
of finite type, where $\frak g$ is simple and $\frak k$ is a
principal $\frak{sl}(2)$-subalgebra.   By Proposition 5, $\frak g[M]$ is a semisimple proper subalgebra of $\frak g$,
and no simple factor of $g[M]$ is isomorphic to $\frak{sl}(2)$.
Theorem 9 (i) is then equivalent to the claim that
$\frak k = \frak g [M]$.  Assume, to the contrary, that the
inclusion $\frak k\subset\frak g [M]$ is proper.  Then $M$ has finite type over $\frak g [M]$, and as $\dim
M = \infty$, infinitely many $\frak g [M]$-types occur in $M$.
Now Proposition 4 implies that the multiplicity of some $\frak
k$-type $W_0$ with $\dim W_0 < b (\frak g [M])$ is infinite in
$M$.  This proves Theorem 9 (i). \vskip .10in

{\bf 4.2. Generalities on cohomological induction.} Let $\frak t
\subset\frak k\subset\frak g$ be a triple of reductive Lie
algebras such that $\frak t$ is reductive in both $\frak k$ and
$\frak g$ and $\frak k$ is reductive in $\frak g$.  Let $C(\frak
g,\frak k)$ (respectively, $C(\frak g,\frak t)$) be the category
of $(\frak g, \frak k)$-modules (resp. $(\frak g,\frak
t)$-modules) which are semisimple as $\frak k$-modules (resp.
$\frak t$-modules).  Note that, by Proposition 1, every
irreducible $(\frak g,\frak k)$-module (resp. $(\frak g,\frak
t)$-module) is an object of $C(\frak g,\frak k)$ (resp. $C(\frak
g,\frak t)$).  Furthermore, it is well known that
$$\align
&\Gamma_{\frak k,\frak t} : C(\frak g, \frak t)\rightsquigarrow C(\frak g,\frak k) \\
&V\rightsquigarrow \Sigma \Sb W\subset V, \dim W =1, \dim U(\frak
k)\cdot W < \infty\endSb ~ W
\endalign
$$
is a well-defined left exact functor.  We denote by $R^i
\Gamma_{\frak k,\frak t}$ its $i$-th right derived functor.  The
following theorems summarize some important properties of the
functors $R^i \Gamma_{\frak k,\frak t}$. \vskip .10in

\proclaim{Theorem 10} \cite{EW}  Let $M$ be an object of $C(\frak
g,\frak t)$. \roster \item"{(i)}"  $R^i\Gamma_{\frak k,\frak
t}(M) = 0$ for $i > \dim \frak k - \dim \frak t$. \item"{(ii)}"
If $M$ has finite type over $\frak t$, then $R^i \Gamma_{\frak
k,\frak t}(M)$ is a $(\frak g,\frak k)$-module of finite type
over $\frak k$ for every $i\geq 0$. \item"{(iii)}"  If $M$ has
finite type over $\frak t$, then for each $i\leq \dim\frak k -
\dim \frak t$, there is a natural isomorphism of $(\frak g,\frak
k)$-modules
$$
R^i \Gamma_{\frak k, \frak t}(M) \cong (R^{\dim\frak k-\dim\frak t
-i} \Gamma_{\frak k,\frak t} (M^*_{\frak t}))^*_{\frak k},
$$
where $N^*$ denotes the $\frak g$-module dual to a $\frak
g$-module $N$, and $N^*_{\frak t}$ (respectively, $N^*_{\frak
k}$) stands for the maximal submodule of $N^*$ which is an object
of $C(\frak g,\frak t)$ (resp., $C(\frak g,\frak k)$).
\endroster
\endproclaim

\proclaim{Theorem 11} \cite{V}  Suppose $M$ is an object of
finite type in $C(\frak g, \frak t)$ and $W$ is an irreducible
finite dimensional $\frak k$-module.  Then
$$\align
&\sum\limits_{i} (-1)^i \dim \text{Hom}_{\frak k} (W, R^i \Gamma_{\frak k,\frak t} (M)) = \\
& \sum\limits_{i} (-1)^i \dim \text{Hom}_{\frak t} (W\otimes
\Lambda^i (\frak k/\frak t), M),
\endalign
$$
where $\Lambda^i$ stands for $i$-th exterior power.
\endproclaim

Finally, assume that $\frak t$ is abelian and let $\frak h$ be a
Cartan subalgebra of $\frak g$ with $\frak h \supset\frak t$. Fix
a Borel subalgebra $\frak b\subset\frak g$ with $\frak
h\subset\frak b$ and let $M(\lambda)$ denote the Verma module
$U(\frak g)\bigotimes\limits_{U(\frak b)} \bold C_\lambda $,
where $\bold C_\mu$ is the $1$-dimensional $\frak h$-module on
which $\frak h$ acts via $\mu:\frak h\to\bold C$. Then, for each
$i\geq 0$, we define the family of $(\frak g,\frak k)$-modules
$A^i(\lambda)$ by setting
$$
A^i(\lambda):= R^i \Gamma_{\frak k,\frak t} (M(\lambda)).
$$
In general, $A^i(\lambda)$ need not be of finite type over $\frak
k$ as $M(\lambda)$ may not be of finite type over $\frak t$.
However, to ensure this latter condition, it suffices to assume
that $\frak t$ contains an element with strictly positive real
eigenvalues in $\frak n := [\frak b,\frak b]$.  Under this assumption, the
construction of $A^i(\lambda)$ is a close analog of a
construction in Harish-Chandra module theory known as
cohomological induction, cf. \cite{KV}. \vskip .20in

{\bf 4.3.  Proof of Theorem 9 (ii)}. Now let $\frak
k\subset\frak g$ be a principal $\frak{sl}(2)$-subalgebra, and
$\{ e,h,f\}$ be a standard basis for $\frak k$.  Then $h$ is a
regular semisimple element in $\frak g$, and $ad_h:\frak
g\to\frak g$ has even integral eigenvalues, since as a $\frak
k$-module $\frak g$ is isomorphic to a direct sum of irreducible
odd dimensional modules, see \cite{D2}.  We set $\rho_{\frak k}:
= \frac{\alpha}{2}$, where $\alpha$ is the root of $\frak k$ with
root space $\langle e\rangle_{\bold C}$.

Let $\frak t := \langle h\rangle_{\bold C}$.  Then $\frak h :=
C(\frak t)$ is a  Cartan subalgebra of $\frak g$ with $\frak
h\supset \frak t$.  Define the Borel subalgebra $\bar\frak b
\supset\frak h$ as the sum of all nonnegative eigenspaces of
$ad_h$, and put $\bar\frak n := [\bar\frak b, \bar\frak b]$. Then
$\bar\frak n = \langle e\rangle_{\bold C} + \bar\frak n\cap\frak
k^\perp$, where $\perp$ stands for orthogonal complement with
respect to the Killing form on $\frak g$. Furthermore, for any
weight $x$ of $\bar\frak n\cap\frak k^\perp$ as a module over
$\langle h\rangle_{\bold C}$, $x(h)$ is a positive even integer.
Finally, let $\frak b$ be the Borel subalgebra opposite to $\frak
b$, i.e. $\frak b\cap\bar\frak b = \frak h$.

Consider the $(\frak g,\frak k)$-modules $A^0(\lambda),
A^1(\lambda)$ and $A^2(\lambda)$ defined as above (by generalized cohomological
induction) for the fixed $\frak t, \frak k, \frak b$ and $\frak
h$.  As the eigenvalues of $ad_{-h}$ on $\frak n = [\frak b,\frak b]$ are positive,
the $\frak g$-modules $A^0(\lambda), A^1(\lambda), A^2(\lambda)$
have finite type over $\frak k$.  We claim first that
$A^0(\lambda) = 0$ for any $\lambda$.  Indeed, by definition
$A^0(\lambda) = \Gamma_{\frak k,\frak t}(M(\lambda))$ and thus
$\frak k, \frak b\subset\frak g [\Gamma_{\frak k,\frak
t}(M(\lambda))]$.  Since $\frak k$ is a stem subalgebra, $\frak
k$ and $\frak b$ generate $\frak g$, i.e. $\frak g [A^0(\lambda)]
= \frak g$.  Hence $A^0(\lambda) = 0$ as $M(\lambda)$ has no
nontrivial integrable $\frak g$-submodules.

We claim next that $A^2(\lambda) = 0$ for any nonintegral
$\lambda$.  This follows from Theorem 10 (iii), once we check
that $\Gamma_{\frak k,\frak t} (M(\lambda)^*_{\frak t}) = 0$.  We
easily see that $M(\lambda)^*_{\frak t} = M(\lambda)^*_{\frak
h}$.  Since $\lambda$ is nonintegral, $M(\lambda)$ has no finite
dimensional nontrivial quotient $\frak g$-module and thus
$M(\lambda)^*_{\frak h}$ has no nontrivial integrable
$\frak g$-submodule.

The following proposition is a statement about the structure of
$A^1(\lambda)$ as a $\frak k$-module and is proved by a
non-difficult computation based on Theorem 11. If $\mu\in\frak
t^*$ is a dominant integral weight for $\frak k$, let $W(\mu)$ be
the irreducible finite dimensional $\frak k$-module with highest
weight $\mu$.  For $\nu\in \frak t^*$, let $P_{\bar\frak
n\cap\frak k^{\perp}} (\nu)$ be the  multiplicity of $\nu$ in the symmetric algebra on 
$\bar\frak n\cap\frak k^\perp$.

\proclaim{Proposition 6}  For any dominant integral weight $\mu$
for $\frak k$ and for any non $\frak g$-integral $\lambda\in\frak
h^*$, we have
$$
\dim \text{Hom}_{\frak k}(W(\mu), A^1(\lambda)) = P_{\bar\frak n\cap\frak
k^{\perp}} (\mu - \lambda\mid_t + 2\rho_{\frak k}) - P_{\bar\frak
n\cap\frak k^{\perp}} (-\mu - \lambda\mid_{\frak t}). \tag3
$$
\endproclaim

Note that if $0 \neq \lambda (h)\in \Z_+$, then the second term on the right hand side of (3) vanishes.  This leads to the following corollary.

\proclaim{Corollary 1}  Suppose $\lambda (h) - 2 = n\in\bold
Z_+$. Then $W(n\rho_{\frak k})$ occurs with multiplicity $1$ in
$A^1(\lambda)$, and $m \geq n$ whenever $W(m\rho_{\frak k})$
occurs in $A^1(\lambda)$.
\endproclaim

The proof of Theorem 9(ii) is now immediate. Let the highest weight of
$Z$ be $n\rho_{\frak k}$ for $n\in\bold Z_+$.  Choose
$\lambda\in\frak h^*$ so that $\lambda$ is not $\frak g$-integral
but such that $\lambda (h) - 2 = n$.  By Corollary 1, $Z\cong
W(n\rho_{\frak k})$ occurs with multiplicity $1$ in
$A^1(\lambda)$ and $m > n$ whenever $W(m\rho_{\frak k})$ occurs
in $A^1(\lambda)$.  Define $X$ as an irreducible quotient of the
$\frak g$-submodule of $A^1(\lambda)$ generated by $Z$.  Then $Z$
occurs with multiplicity $1$ in $X$ and $m > n$ whenever
$W(m\rho_{\frak k})$ occurs in $X$.

As $\lambda$ is not $\frak g$-integral, the infinitesimal
parameter of $M(\lambda)$, and hence (see \cite{V}) of
$A^1(\lambda)$ and $X$, is not integral.  Thus $X$ is infinite
dimensional. \qed \vskip .15in

{\bf Remark.}  If we fix $n$ and set $\lambda (h) - 2 = n$, then
$\lambda$ still has $\ell -1$ complex parameters $(\ell =
\dim\frak h)$. It will be quite interesting to determine the
$\frak g$-module structure of $A^1(\lambda)$ as $\lambda$ varies
continuously with $\lambda(h)$ fixed.  Will $A^1(\lambda)$ have
finite length for all $\lambda$ ?  This is known (see \cite{V})
to be true if $\frak g = \frak{sl}(3)$.

\vskip .25in

\noindent Department of Mathematics

\noindent Yale University

\noindent New Haven, CT 06520

\noindent E-mail Address: penkov\@math.yale.edu

\noindent {\bf Permanent Address:}

\noindent Department of Mathematics

\noindent University of California at Riverside

\noindent Riverside, CA  92521
\vskip .20in

\noindent Department of Mathematics

\noindent Yale University

\noindent New Haven, CT  06520

\noindent Email Address:  gregg\@math.yale.edu

\end